\documentclass[12pt]{amsart}

\usepackage{amscd,amssymb}

\def\draft{y}
\def\printname#1{
	\if\draft y
		\smash{\makebox[0pt]{\hspace{-0.5in}
			\raisebox{8pt}{\tt\tiny #1}}}
	\fi
}
\catcode`\@=11
\long\def\@makecaption#1#2{%
    \vskip 10pt
    \setbox\@tempboxa\hbox{
      \small\sf{\bfcaptionfont #1. }\ignorespaces #2}%
    \ifdim \wd\@tempboxa >\captionwidth {%
        \rightskip=\@captionmargin\leftskip=\@captionmargin
        \unhbox\@tempboxa\par}%
      \else
        \hbox to\hsize{\hfil\box\@tempboxa\hfil}%
    \fi}
\font\bfcaptionfont=cmssbx10 scaled \magstephalf
\newdimen\@captionmargin\@captionmargin=2\parindent
\newdimen\captionwidth\captionwidth=\hsize
\catcode`\@=12

\newtheorem{theorem}{Theorem}[section]

\theoremstyle{definition}

\theoremstyle{remark}

\numberwithin{equation}{section}

\newcommand{\Z}{{\mathbb Z}}

\newcommand{\R}{{\mathbb R}}
\newcommand{\C}{{\mathbb C}}

\newcommand{\OO}{{\mathcal {O}}}
\newcommand{\M}{{\mathcal {M}}}
\newcommand{\T}{{\mathcal {T}}}
\newcommand{\m}{{\mathfrak {m}}}

\DeclareMathOperator{\sign}{{sign}}
\DeclareMathOperator{\End}{{End}}
\DeclareMathOperator{\Eul}{{Eul}}
\DeclareMathOperator{\rk}{{rk}}
\DeclareMathOperator{\mf}{{Mon_F}}

\begin{document}

\title[Absolute torsion and eta-invariant] %
{Absolute torsion and eta-invariant}

\author[M.~Farber]{M.~Farber}
\address{Department of Mathematics, Tel Aviv University, Tel Aviv, 69978, Israel}
\email{Farber@math.tau.AC.IL}
\date{\today}

\thanks{Partially supported by the US - Israel Binational Science Foundation and by the Minkowski Center for Geometry}


\maketitle

In a recent joint work with V. Turaev \cite{FT1}, 
we defined a new concept of combinatorial 
torsion which we called {\it absolute torsion}. Compared with the
classical Reidemeister torsion, it has the advantage of having a well-determined sign. Also, the absolute torsion is defined for arbitrary orientable flat vector bundles, and not only for unimodular ones, 
as is classical Reidemeister torsion.

In this paper I show that the sign behavior of the absolute torsion, under a continuous deformation of the flat bundle, is determined by the eta-invariant and the Pontrjagin classes.

\section{\bf  A review of absolute torsion}
\label{s:abs}

In this section we briefly review the main properties of {\it absolute torsion}, which was introduced in \cite{FT1}.

Let $X$ be a closed {\it oriented}
PL manifolds of odd dimension $m=\dim X$, and let $F\to X$ be a flat complex vector bundle over 
$X$. Suppose that the following conditions hold:

{\it (i) The Stiefel-Whitney class $w_{m-1}(X)\in
H^{m-1}(X,\Z_2)$ vanishes;

(ii) The first Stiefel-Whitney class $w_1(F)$, viewed as a homomorphism $H_1(X;\Z)$ $\to \Z_2$,
vanishes on the 2-torsion subgroup of  $H_1(X;\Z)$.}

Note that condition (i) is automatically satisfied in the case 
when $m\equiv 3\mod 4$,
as proven by W. Massey \cite{Ma}. Condition (ii) holds for any orientable bundle $F$.
(ii) also holds for any $F$, assuming that $H_1(X)$ has no 2-torsion.

Under the above conditions we defined in \cite{FT1} an element of the determinant line of homology
\begin{eqnarray}
\T(F)\in \det\, H_\ast(X;F),\label{abstor}
\end{eqnarray}
which we called {\it absolute torsion.} The term intends to emphasize that this invariant is independent
of additional choices, such as an Euler structure.

The construction of absolute torsion is based on the construction of torsion of Euler 
structures, which was initiated by V. Turaev in \cite{T1}, \cite{T2} 
and studied further in \cite{FT} in the context of flat vector bundles and determinant lines.
In general, the construction of torsion of Euler structures allows to control the indeterminacy 
contained in the classical
construction of the Reidemeister torsion. However, the resulting invariant really depends on 
the choice of Euler structure and also on the 
choice of homological orientation. In \cite{FT1} we showed that under
condition {\it (i)} above there exist canonical choices for the Euler structure, which are determined by the requirement that their characteristic classes vanish. 
In fact, there exist only finitely many {\it canonical Euler structures}; they 
are parametrized by the 2-torsion subgroup of $H_1(X)$.  Under condition {\it (ii)} the torsion
corresponding to any of the canonical Euler structures will be the same.
There is also a canonical way of producing a homological orientation using the orientation of $X$.

We proved in \cite{FT1} that the absolute torsion (\ref{abstor}) is well-defined, is combinatorially
invariant, has no sign indeterminacy, but may depend on the orientation of $X$. More precisely, reversing the orientation of $X$
multiplies the absolute torsion (\ref{abstor}) by the number
$(-1)^{\rk(F)\cdot s\chi(X)},$
where $s\chi(X)$ denotes the {\it semi-characteristic} of $X$, i.e.
\[s\chi(X) = \sum_{i=0}^{(m-1)/2} b_{2i}(X).\]
Hence the absolute torsion is independent of the orientation of $X$ if one of the numbers $\rk(F)$ and $s\chi(X)$ is even.

Another important property of absolute torsion is {\it reality}, which was established in \cite{FT1},
Theorem 3.8. It claims that if $F$ is a flat complex vector bundle 
over a closed odd dimensional manifold $X$, such that conditions {\it (i) and (ii)} are satisfied,
and $F$ admits a flat Hermitian metric, then the absolute torsion $\T(F)\in \det\, H_\ast(X;F)$
is real. This means that 
\[\overline {\T(F)} = \T(F),\]
where the bar denotes a canonical involution on the
determinant line $\det\, H_\ast(X;F)$, cf. \cite{FT1}, \S 2. 

Note that assuming that the flat bundle $F$
is acyclic, the determinant line $\det\, H_\ast(X;F)$ identifies canonically with the field of complex numbers
\[\det\, H_\ast(X;F)\, \simeq \, \C,\]
and the involution
on the determinant line turns into the usual complex conjugation, cf. \cite{FT1}, Lemma 2.2. Hence
in the acyclic case the absolute torsion $\T(F)$ is a real valued function of a flat Hermitian vector bundle $F$, and our purpose in this paper is to describe its sign. We will see that this description involves the eta-invariant \cite{APS} and also
the Pontrjagin classes. Note that the absolute value of the absolute torsion in this situation coincides
with the Ray-Singer analytic torsion, as follows from Theorem 10.2 of \cite{FT}.

Note, that the absolute torsion can be viewed as a high dimensional generalization
of the Conway polynomial. In fact, we proved in \cite{FT1} that the function $F\mapsto \T(F)$
determines the Conway polynomial of a knot in a 
3-sphere, where $F$ is a flat line bundle over the
closed 3-manifold obtained as the 
result of 0-framed surgery along the knot.

\section{\bf Main Theorem in case $\dim X\equiv 3\mod 4$}
\label{eta1}

\subsection{Assumptions}
Let $X$ be a closed oriented smooth manifold of odd dimension $m=\dim X$. In this subsection we will assume that $m\equiv 3\mod 4$. By a theorem of W. Massey \cite{Ma}, the Stiefel - Whitney class $w_{m-1}(X)$ vanishes. 
Let $F$ be Hermitian vector bundle $F$ over $X$, such that the Stiefel - Whitney class $w_1(F)$ is trivial on the 2-torsion subgroup of $H_1(X)$.

Suppose that for $t\in (a, b)$ we are given a {\it real analytic} curve of flat Hermitian connections $\nabla_t$ on $F$. The analyticity of $\nabla_t$ with respect to $t$
is understood as follows. Fix $t_0\in (a,b)$; then we have $\nabla_t = \nabla_{t_0} + \Omega_t$, where $\Omega_t\in A^1(X;\End(F))$ is a 1-form with
values in the bundle of endomorphisms of $F$. The curve of connections $\nabla_t$ 
is analytic if the curve
\[(a,b)\to A^1(X;\End(F)),\quad t\mapsto \Omega_t\]
is analytic with respect to any Sobolev
norm. We refer to \cite{F}, sections 2.3, 2.4 and \cite{FL}, sections 3.3, 3.4 for detailed definitions 
of these notions. 

Let $F_t$ denote the flat Hermitian bundle $(F, \nabla_t)$, $t \in (a,b)$. 
We will assume that the homology $H_\ast(X;F_t)$ vanishes for a generic 
$t\in (a,b)$, i.e. for all $t\in (a,b)-S$, where $S\subset (a,b)$ is a finite subset. Then the absolute torsion $t\mapsto \T(F_t)$ gives a real valued function of $t\in (a,b)$, $t\notin S$. 
The function $\T(F_t)$ is real analytic and is nonzero on $(a,b)-S$ 
as follows directly from the definitions of \cite{FT1}. 

Fix a Riemannian metric on $X$.
The flat Hermitian connection $\nabla_t$, where $t\in (a,b)$, determines
its eta-invariant $\eta_t=\eta(F_t)$, cf. \cite{APS}. Recall that
the eta-invariant is defined as the value of the
eta-function $\eta(s)$ at the origin $s=0$. 
The eta-function $\eta(s)$ is obtained by an analytic continuation
of 
$$\eta(s) = \sum_{\lambda\ne 0} \sign(\lambda)|\lambda|^{-s}, \quad \Re(s) \quad \text{large},$$
where $\lambda$ runs over the eigenvalues of the self-adjoint elliptic operator 
\[B: A^{ev}(X;F)\to A^{ev}(X;F)\]
acting on smooth forms
with values in $F$ of even degree given by
$$B\phi = i^{r}(-1)^{p+1}(\ast\nabla - \nabla\ast)\phi, \quad \phi\in A^{2p}(X;F).$$
Here $m=\dim X = 2r-1$ and the star denotes the Hodge duality operator. 

The following is the Main result of the paper in the case $m\equiv 3\mod 4$.

\begin{theorem}
\label{theorem1}
Let $X$ be a closed oriented Riemannian manifold of dimension 
$m\equiv 3\mod 4$ and let $\nabla_t$ 
be an analytic family of flat Hermitian connections, $t\in  (a,b)$,
on an orientable vector bundle $F$ over $X$. Assume that the flat bundle $(F, \nabla_t)$ is acyclic
for all $t\in (a,b)$, $t\notin S$, where $S\subset (a,b)$ is a finite subset.
Then the complex number
\begin{eqnarray}
\sign(\T(F_t)) \exp({i\pi\eta(F_t)/2)}\in \C\label{ray1}
\end{eqnarray}
is independent of $t\in  (a,b)$, $t\notin S$.
\end{theorem}

Here $\sign(\T(F_t))$ denotes the sign of the absolute torsion $\T(F_t)$.

Proof is given in section \ref{proof} below.

\section{\bf Main Theorem in case $\dim X\equiv 1\mod 4$}
\label{eta2}

The difference between the cases $m\equiv 3\mod 4$ and $m\equiv 1\mod 4$ consists in the following.
In the case $m\equiv 3\mod 4$ the eta-invariant $\eta(F_t)$ is locally constant and may have only integral jumps
at the points $t\in (a,b)$, where the acyclicity is violated.
On the contrary, in the case $m\equiv 1\mod 4$ the eta-invariant $\eta(F_t)$ as a function 
of $t$ behaves ``linearly'' between the jump points. This explains the difference with the case $m\equiv 3\mod 4$. 

In order to state the Main result in the case $m\equiv 1\mod 4$ we need the following
construction.

\subsection{The argument class} 
\label{arg}
Let $F$ be a flat Hermitian bundle over $X$. For any closed curve
$\gamma$ in $X$, starting at the base point, we have the operator of parallel transport
$\mf(\gamma) :F_0\to F_0$ along $\gamma$,
where $F_0$ denotes the fiber over the base point. The linear map $\mf(\gamma)$ is unitary, and 
hence its determinant lie on the unit circle.
We will denote 
by $\mathbf {Arg}_F$ a cohomology class $\mathbf{Arg}_F\in H^1(X;\R/\Z)$ 
with the following property: for any closed curve
$\gamma$ in $X$ 
holds
\[\det (\mf(\gamma)) \, = \, \exp( 2\pi i\cdot \langle \mathbf{Arg}_F,[\gamma]\rangle).\]

\begin{theorem}
\label{theorem2} 
Let $X$ be a closed oriented Riemannian manifold of dimension $m\equiv 1\mod 4$. Suppose that 
$H_1(X)$ has no 2-torsion and the Stiefel-Whitney class $w_{m-1}(X)\in H^{m-1}(X;\Z_2)$ 
vanishes.
Let $\nabla_t$ 
be an analytic family of flat Hermitian connections, $t\in  (a,b)$,
on a vector bundle $F$ over $X$, such that $H^\ast(F, \nabla_t)=0$ for all $t\in (a,b)-S$, where $S\subset (a,b)$ is a finite subset. Then 
\begin{enumerate}
\item[(i)] the quantity
\begin{eqnarray}
\exp(i\pi \langle\mathbf {Arg}_{F_t}\cup L(X),[X]\rangle)\quad \in S^1\subset \C\label{argument}
\end{eqnarray}
is well defined (cf. \S 4.4 below) and
represents a continuous function of the parameter $t\in (a,b)$;
\item[(ii)] the complex number
\begin{eqnarray}
\sign(\T(F_t))\exp({i\pi\eta(F_t)/2)}
\exp(i\pi \langle\mathbf {Arg}_{F_t}\cup L(X),[X]\rangle)\, \in \C\label{ray2}
\end{eqnarray}
is independent of $t\in  (a,b),  t\notin S$.
\end{enumerate}
\end{theorem}

In formula (\ref{ray2}) $L(X)$ denotes the Hirzebruch polynomial in the Pontrjagin classes. 

The numerical value of the eta-invariant  $\eta_t$ depends on the 
choice of the Riemannian metric on $X$. Therefore the arguments of the rays in (\ref{ray1}) and (\ref{ray2}) will
depend on the choice of the metric on $X$. Also, the orientation of $X$ used in order to define
$\eta(F_t)$ and also (\ref{argument}), so the arguments of the rays in  (\ref{ray1}) and  (\ref{ray2}) will depend also on the choice of the orientation of $X$. 

A proof of Theorem \ref{theorem2} will be given in the next section.

There exist a version of Theorem \ref{theorem2}, which is similar to Theorem \ref{theorem1}. Here
we do not assume absence of 2-torsion in $H_1(X)$, although we require that the family of flat connections $\nabla_t$ have trivial determinant of the monodromy $\det \mf$ (in other words, each 
$\nabla_t$ is a $SU$-connection).

\begin{theorem}
\label{theorem3}
Let $X$ be a closed oriented Riemannian manifold of dimension $m\equiv 1\mod 4$ such that 
the Stiefel-Whitney class $w_{m-1}(X)$ vanishes. Let $\nabla_t$ 
be an analytic family of flat $SU$-connections, $t\in  (a,b)$,
on an orientable vector bundle $F$ over $X$. Assume that the flat bundle $(F, \nabla_t)$ is acyclic
for all $t\in (a,b)$, $t\notin S$, where $S\subset (a,b)$ is a finite subset.
Then the complex number
\begin{eqnarray*}
\sign(\T(F_t))\exp({i\pi\eta(F_t)/2)}\end{eqnarray*}
is independent of $t\in  (a,b), t\notin S$.
\end{theorem}

\section{\bf Proofs of Theorems \ref{theorem1}, \ref{theorem2} and \ref{theorem3}} 
\label{proof}

The proofs will be based on the results of \cite{F}, \cite{FL1}, \cite{FL}, where
the behavior of the eta-invariant and the analytic torsion 
under the analytic perturbations of flat connection was described.

\subsection{Signature invariants of deformations}
We will first give a brief summary of the results of \cite{FL1}, \cite{FL}.

There are results of two types in \cite{FL}: the description of the variation of the eta-invariant 
$\eta_t$ modulo $\Z$ (cf. Theorem 7.1 in \cite{FL}, and also the description of the integral jumps
(cf. Theorem 1.5 in \cite{FL}). 

In the case $m\equiv 3\, \mod \, 4$ the eta-invariant $\eta_t$ as a function of $t$ is constant
modulo 1 (may have only integral jumps). The jumps occur only at points where the homology changes; such points are isolated (not generic). If $t=t_0\in S$ is a jump point then the limits 
$\lim_{t\to +t_0}\eta_t= \eta_+$ and $\lim_{t\to -t_0}\eta_t= \eta_-$ exist. The jump across
$\eta_+-\eta_-$ is always even. Theorem 1.5 of \cite{FL} gives a precise formulae for the jump
across $\eta_+-\eta_-$ in homological terms.
To express these formulae one needs the signature invariants $\sigma_1, \sigma_2, \dots$, 
which were constructed in \cite{FL}, \S 2; they are determined in pure homological 
(finite dimensional) terms as invariants of some {\it linking
form}, naturally determined by {\it the deformation of the monodromy representation}, cf. \cite{FL}.
The jump formulae claim that 
\begin{eqnarray}
\eta_+-\eta_-\, =\, 2\sum_{i\,  \text{odd}}\sigma_i\, =\, 2\sigma_{\text{odd}},\label{jump}
\end{eqnarray}
and also, the value of the eta-invariant at the jump point itself $\eta_{t_0}$ is given by
\begin{eqnarray}
\frac{1}{2}(\eta_++\eta_-) -\eta_{t_0} =\sum_{i\, \text{even}} \sigma_i\, =\, \sigma_{\text{even}}.
\label{mid}
\end{eqnarray}
It follows that the values of the 
function $t\mapsto \exp(i\pi\eta_t/2)$ for a generic $t$ always lie on a fixed straight line through
the origin and it changes the phase (the half line), while passing across a jump point $t=t_0$, if and only if the signature $\sigma_{\text{odd}}$ is odd.  

Note that the absolute torsion is continuous, real and nonzero in the intervals between the jump points (where there are no homology changes) and so the sign 
of $\T(F_t)$ is constant between the jump points.
Hence, 
Theorem \ref{theorem1} is equivalent to the inequality
\begin{eqnarray}
\frac{\T(F_{t_0+\delta})}{\T(F_{t_0-\delta})}\cdot\exp(i\pi\sigma_{\text{odd}})\,  > \, 0\label{goal}
\end{eqnarray}
for $\delta>0$ small enough. 

\subsection{Singularity of the torsion} Here we will use the theorem of \cite{F}, which describes
the singularity of the torsion.

Let $t_0\in S \subset (a, b)$ be a fixed point, 
where the homology of flat bundle $H^\ast(X;F_t)$ is nontrivial.
It is clear that as a function of $t$ the absolute torsion $\T(F_t)$ has the form
\begin{eqnarray}
\T(F_t) = (t-t_0)^{\nu}\cdot f(t), \quad \nu \in \Z,\label{exp}
\end{eqnarray}
where $f(t)$ is real analytic and nonzero in a neighborhood of $t_0$; one obtains (\ref{exp}) directly 
from the definition of torsion in \cite{FT}, \cite{FT1}. The exponent $\nu$
describes the singularity of the torsion at $t=t_0$ (the order of zero or pole). For our purposes in this paper it is enough to know the exponent $\nu$ only modulo 2, since the absolute torsion $\T(F_t)$
changes sign while $t$ passes $t=t_0$ if and only if $\nu$ is odd.

It follows that in order to prove (\ref{goal}) we need to show that
\begin{eqnarray}
\nu \equiv \sigma_{\text{odd}}\mod 2\label{goal1}
\end{eqnarray}

We will prove (\ref{goal1}) using Theorem 10.2 of \cite{FT} and Theorem 5.3 of \cite{F}. Since Theorem 10.2 of \cite{FT} operates with
the cohomological torsion and real vector bundles, we need to adjust our notations. Let 
\[\tau^\bullet_\R(F_t) = \tau^{\bullet}(X,\xi;F_t^{\R})\in \det \, 
H^\ast(X;F_t^\R)\simeq \R,\quad t\ne t_0,\]
be the cohomological torsion (cf. \cite{FT}, \S 9) of flat vector bundle $F_t$ considered as a real vector bundle. 
Here $\xi\in \Eul(X)$ denotes a {\it canonical Euler structure},
cf. \cite{FT1}, \S 3; we also assume that $X$ is supplied with {\it the canonical homological orientation}, cf. \cite{FT}, \S3. 

Let 
\[\tau^\bullet_\C(F_t) = \tau^{\bullet}(X,\xi;F_t^{\C})\in \det \, H^\ast(X;F_t)\simeq \C,\quad t\ne t_0,\]
be the (complex) cohomological torsion of $F_t$ (cf. \cite{FT}, \S 9).
By formula (9-3) in \cite{FT} we have
\[|\tau^\bullet_\C(F_t)| = |\T(F_t)|^{-1} = |t-t_0|^{-\nu}\cdot |f(t)|^{-1}.\]

By Theorem 10.2 of \cite{FT} (cf. also the equation (10-3) in \cite{FT}), 
we obtain that for $t\ne t_0$, the analytic torsion $\rho(F_t^\R)$ of the flat bundle $F_t^\R$ equals 
\begin{eqnarray}
\rho(F_t^\R) = |\tau^\bullet_\R(F_t)|.\label{anal}
\end{eqnarray}
Since clearly
\[\tau^\bullet_\R(F_t) = |\tau^\bullet_\C(F_t)|^2, \quad \rho(F_t^\R) = \rho(F_t)^2,\]
we obtain that 
\[\rho(F_t) = |\tau^\bullet_\C(F_t)| = |t-t_0|^{-\nu}\cdot |f(t)|^{-1}, \quad f(t_0)\ne 0.\]

Now we may apply Theorem 5.3 of \cite{F} in order to compute the exponent $\nu$. According to this Theorem, $\nu$ equals the Euler number of the deformation 
\begin{eqnarray}
\chi = \sum_{i=0}^m (-1)^i\dim_\C \mathfrak T^i,\label{chi}
\end{eqnarray}
cf. \cite{F}.  

Recall the notation used in the last formula.
Let $\OO=\OO_{t_0}$ denote the ring of germs at $t=t_0$ of real analytic curves $f: (t_0-\epsilon, t_0+\epsilon)\to \C$. Addition 
and multiplication are given by pointwise operations. $\OO$ is a 
discrete valuation ring; its maximal ideal $\m\subset \OO$ coincides with
the set of all functions vanishing at $t=t_0$. The generator of the 
maximal ideal is the germ of the function $f(t)=t-t_0$.

Given an analytic family of flat connections $\nabla_t$ (where $t\in (a, b)$) on a Hermitian vector bundle $F$, 
one constructs (following \cite{F}, \cite{FL})
a single local system $\OO F$ of free $\OO$-modules of rank $\dim F$ over $X$. To describe it, we denote by $\Pi_t(\gamma): F_p\to F_p$ the 
parallel transport in $F$ with respect to the flat connection $\nabla_t$ along the loop 
$\gamma: [0,1]\to X$, where $p=\gamma(0)=\gamma(1)$ and $F_p$ denotes the fiber above $p$. 
The bundle $\OO F$
as a set consists of germs of analytic curves $g: (t_0-\epsilon', t_0+\epsilon')\to F,$ which are vertical
(i.e. belong to a single fiber). The parallel transport of the bundle $\OO F$ along a closed loop
$\gamma: [0,1]\to X$, where $p=\gamma(0)=\gamma(1)$ is defined as 
$\Pi(\gamma): \OO F_p\to \OO F_p$, where for $g\in \OO F_p$
\[ \Pi(\gamma)(g)(t) = \Pi_t(\gamma)(g(t)),\]
where $t\in (t_0-\epsilon', t_0+\epsilon').$

Consider the cohomology $H^\ast(X;\OO F)$ with coefficients in the local system $\OO F$. Since
the ring $\OO$ is a discrete valuation ring, any cohomology module $H^i(X;\OO F)$ is a direct sum
of its torsion submodule and a free module. Our assumptions (acyclicity of the flat bundle $F_t$ for
generic $t$) imply that the free part of the homology of $H^i(X;\OO F)$ is trivial (cf. Theorem 2.8 
in \cite{F}). We will denote by ${\frak T^i}$ the torsion submodule of 
$H^i(X;\OO F)$. Note that each ${\frak T^i}$ is finite dimensional as a $\C$-vector space.

\subsection{End the proof of Theorem \ref{theorem1}}
Let $\M$ denote the field of fractions of $\OO$. Elements $f\in \M$ could be viewed as germs at
$t=t_0$ of meromorphic curves $(t-t_0)^\nu\cdot f(t)$, where $\nu \in \Z$, $f\in \OO$, $f(t_0)\ne 0$.

In \cite{FL}, \S 1, {\it a linking form} on the middle dimensional torsion submodule
\begin{eqnarray}
{\frak T}^r \otimes {\frak T}^r\, \to \M/\OO,\qquad m=2r-1\label{link}
\end{eqnarray}
was constructed. The signatures $\sigma_1, \sigma_2, \dots$ (which appear in (\ref{jump}) and (\ref{mid}))
are derived from the linking form (\ref{link}). 
From the construction of the signatures $\sigma_i$ in \S 2
of \cite{FL} we easily observe
\begin{eqnarray}
\dim_\C {\frak T}^r \equiv \sigma_{\text{odd}}\mod 2.\label{odd}
\end{eqnarray}

By the Poincar\'e duality we have $\dim_\C \frak T^i = \dim_\C\frak T^{2r-i}$ (cf. \cite{F}, Prop. 7.6).
Hence we obtain combining (\ref{chi}) and (\ref{odd})
\[\nu = \chi \equiv \dim \frak T^r \equiv \sigma_{\text{odd}}\mod 2,\]
which finishes the proof. \qed

\subsection{Proof of Theorem \ref{theorem2}} Let us prove the first statement.
We claim that the quantity
\begin{eqnarray}
\langle{\mathbf{Arg}}_{F_t}\cup L(X),[X]\rangle\label{integ}
\end{eqnarray}
defines a real number modulo $2\Z$. To show this we observe that
the part of the Hirzebruch form $L(X)$ of degree $m-1$ represents an {\it even} cohomology class
(i.e. a class in $2\cdot H^{m-1}(X;\Z)$). This follows from our assumptions that the 
Stiefel - Whitney class $w_{m-1}(X)$ vanishes and also that $H_1(X)$ has no 2-torsion. 
Indeed, any cycle $z\in H_{m-1}(X)$ can be represented by an oriented codimension one submanifold $M_z\subset X$ and the value 
$\langle L(X),z\rangle$ equals the signature of $M_z$. To show that the signature of $M_z$ is always even, it is enough to show that the Euler characteristic of $M_z$ is even. 
But
$\chi(M_z) \equiv  \langle w_{m-1}(X), z\rangle \equiv 0\mod 2,$
which is a part of our data.

This shows that the class $1/2\cdot L^{m-1}(X)\in H^{m-1}(X;\Z)$ exists and is unique, since we assume that $H^{m-1}(X;\Z)\simeq H_1(X)$ has no 2-torsion.

Therefore, we obtain that the quantity 
\[\exp(\pi i \langle \mathbf{Arg}_F\cup L(X),[X]\rangle)=
\exp(2\pi i \langle \mathbf{Arg}_F\cup (1/2\cdot L(X)),[X]\rangle)\]
is well defined, since the cup-product $\mathbf{Arg}_F\cup (1/2\cdot L(X))$ is a well defined class in
$H^m(X;\R/\Z)$. The fact that (\ref{argument}) depends continuously on $t$ is now obvious. This proves the first statement of Theorem \ref{theorem2}.

Consider the quantity
\begin{eqnarray}
\hat\eta_t\, =\, \eta_t + 2\langle\mathbf{Arg}_{F_t}\cup L(X),[X]\rangle \in \R/4\Z.\label{hat}
\end{eqnarray}
Theorem 7.1 of \cite{FL} claims that $\hat\eta_t$, as a function of $t\in (a, b)$, is locally constant.
Hence $\hat\eta_t$ may have only jumps when parameter $t$ passes over 
the finite subset $S\subset (a,b)$. 

Let $t_0\in S$ be one of the jump points. Define $\hat\eta_+$ as the value of $\hat\eta_{t_0+\delta}$
and define $\hat\eta_-$ as $\hat\eta_{t_0-\delta}$, where $\delta>0$ is small enough. Then we have the
following jump formular:
\begin{eqnarray}
\hat\eta_+ - \hat\eta_- = 2\sigma_{\text{odd}}\mod 4.\label{jump2}
\end{eqnarray}
Indeed, as we showed above the quantity (\ref{integ}) is determined in $\R/2\Z$, and hence the
second summand in (\ref{jump2}) has indeterminacy in $4\Z$. Now (\ref{jump2}) follows from 
(\ref{jump}). 

In order to show that 
$\T(F_t)\cdot \exp(\pi i\hat\eta_t/2)$
always stays on the same ray from the origin, it is enough (because of (\ref{jump2})) to check the inequality (\ref{goal}).

However the proof of (\ref{goal}), given in the case $m\equiv 3\mod 4$, works as well 
(with no changes)
in the case $m\equiv 1\mod 4$. \qed

\subsection{Proof of Theorem \ref{theorem3}} It is identical to the proof of Theorem \ref{theorem1}
if one observes that the eta-invariant $\eta_t$, considered modulo $\Z$, is constant along paths of 
$SU$-connections (as follows from Theorem 7.1 in \cite{FL}); hence $\eta_t$ is constant between the
jump points.

\bibliographystyle{amsalpha}

\enddocument